\documentstyle[12pt,twoside, epsf]{article}
 \newtheorem{th}{Theorem}

 \newtheorem{defn}{Definition}

 \newtheorem{rem}{Remark}

\makeatletter
\long\def\M#1{\leavevmode\setbox\@tempboxa\hbox{#1}\@tempdima\fboxrule
     \advance\@tempdima \fboxsep \advance\@tempdima \dp\@tempboxa
    \hbox{\lower \@tempdima\hbox
   {\vbox{\hrule \@height \fboxrule
           \hbox{  \hskip\fboxsep
           \vbox{\vskip\fboxsep \box\@tempboxa\vskip\fboxsep}\hskip
                  \fboxsep\vrule \@width \fboxrule}%
                   }}}}

\makeatother

\def\picill#1by#2(#3)
{\vbox to #2
{\hrule width #1 height 0pt depth 0pt
\vfill\epsffile{#3}}}

\let \ttorg \tt \def \tt{\ttorg \obeyspaces}

\begin{document}

\pagestyle{myheadings}

\title{\bf The $L$--Move and Virtual Braids}

\author{Louis H. Kauffman\\
Department of Mathematics,\\
University of Illinois at Chicago, \\
851 South Morgan Street.\\
Chicago, Illinois 60607-7045,USA \\
E-mail: kauffman@uic.edu\\and\\
Sofia Lambropoulou\\
Department of Mathematics, \\ 
National Technical University of Athens, \\
Zografou campus, GR-157 80 Athens, Greece \\
E-mail: sofia@math.ntua.gr }
 
\maketitle
  
\thispagestyle{empty}

 \section{Introduction}

In this paper we sketch our proof \cite{VL} of a Markov Theorem for the virtual braid group.
This theorem gives a result for virtual knot theory that is analogous to the result of the Markov Theorem
in classical knot theory. We have that every virtual link is isotopic to the closure of a virtual
braid, and that two virtual links, seen as the closures of two virtual braids, are isotopic if and only if
the braids are related by a set of moves. These moves are described in the paper.

In this paper we shall follow the ``$L$--Move" approach to the Markov Theorem (see \cite{Ka,Ka1} for a different
approach).  An $L$--move is a very simple uniform move that can be applied anywhere in a braid to produce a braid with
the isotopic closure. It consists in cutting a strand of the braid and taking the top of the cut to the bottom of the
braid (entirely above or  entirely below the braid) and taking the bottom of the cut to the top of the braid
(uniformly above or below in correspondence with the choice for the other end of the cut). One then proves
that two virtual braids have isotopic closures if  and only if they are related by a series of $L$--moves. Once this
$L$--Move Theorem is established, we can reformulate the result in various ways, including a more algebraic
Markov Theorem that uses conjugation and stabilization moves to relate braids with isotopic
closures. This same approach can be applied to other categories such as welded braids and flat virtual braids
(see \cite{VL}).
\bigbreak

 We first give a quick sketch of virtual knot theory, and then state our Markov Theorem and the definitions that 
support it. The reader interested in seeing the details of this approach should consult \cite{KL} and
\cite{VL}.
\bigbreak

\section{Virtual Knot Theory}

Virtual knot theory is an extension of classical diagrammatic knot theory. In this extension one 
adds  a {\em virtual crossing} (see Figures $1$ and $4$) that is neither an over-crossing
nor an under-crossing.  A virtual crossing is represented by two crossing arcs with a small circle
placed around the crossing point.  
\smallbreak

Virtual diagrams can be regarded as representatives for 
oriented Gauss codes (Gauss diagrams) \cite{VKT,GPV}. Some Gauss codes have planar
realizations, and these correspond to classical knot diagrams. Some codes do not have planar
realizations. An attempt to embed such a code in the plane leads to the production of the virtual
crossings. 
\smallbreak

Virtual knot theory can be interpreted as embeddings of  links in
thickened surfaces, taken up to addition and subtraction of empty handles. 
In this way, we see that this theory is a natural chapter in three-dimensional topology.
\bigbreak

Isotopy moves on virtual diagrams generalize the ordinary Reidemeister moves for classical knot and
link diagrams.  See Figure 1, where all variants of the moves should be considered. 
\bigbreak

$$ \picill5inby3.6in(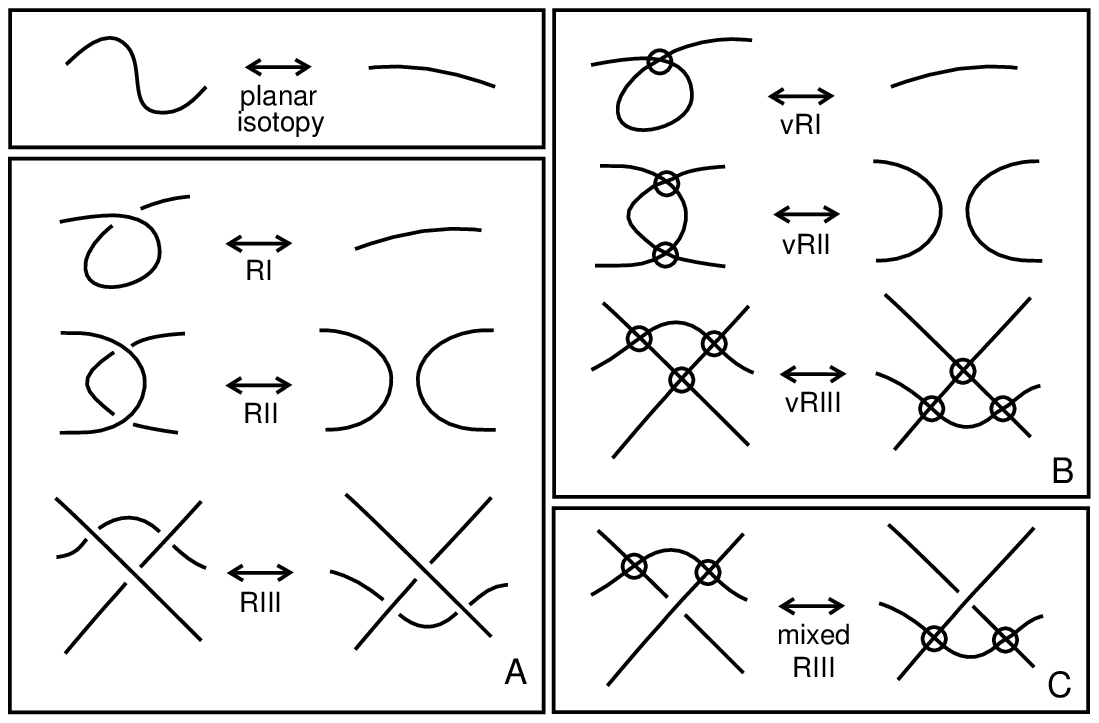)  $$

\begin{center} 
{ \bf Figure 1 --  Reidemeister Moves for Virtuals} 
\end{center}

$$ \picill3inby1.2in(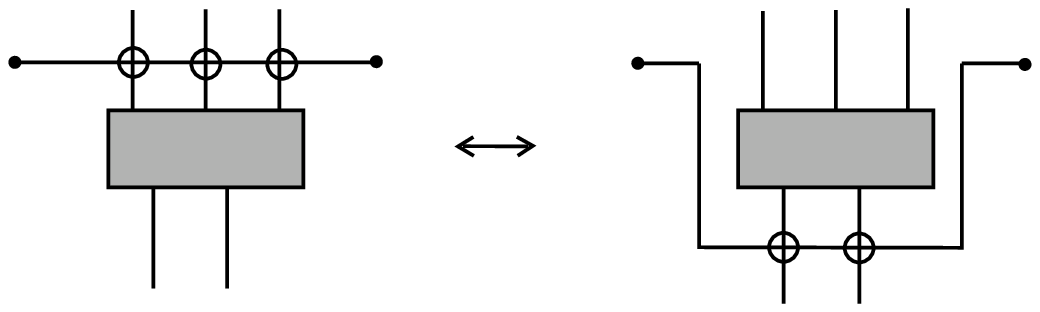)  $$

\begin{center} 
{ \bf Figure 2 -- The Detour Move} 
\end{center}

Equivalently, virtual isotopy is generated by classical Reidemeister moves
and the {\it detour move} shown in Figure 2.  
\bigbreak
 
$$ \picill3inby1in(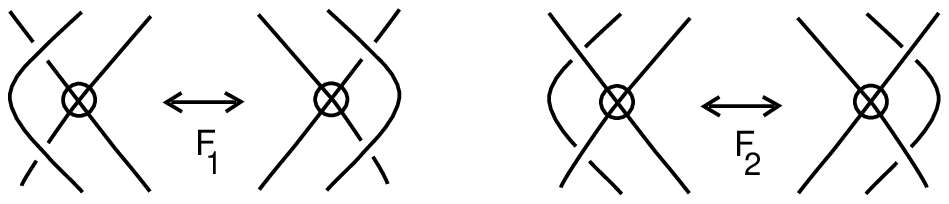)  $$

\begin{center} 
{ \bf Figure 3 -- The Forbidden Moves} 
\end{center}

The moves shown in Figure 3 are forbidden in virtual knot theory (although not in some of its variants such
as welded links). The forbidden moves are not consequences of the notion of virtual equivalence. In working with
the Markov Theorem for virtual knots and links, we have to respect these constraints.
\bigbreak
 
We know \cite{VKT,GPV} that classical knot theory embeds faithfully in virtual knot theory. That is, if two 
classical knots are isotopic through moves using virtual crossings, then they are isotopic as
classical knots via standard  Reidemeister moves.  With this approach, one can generalize many 
structures in classical knot theory to the virtual domain, and use the virtual knots to test the
limits of classical problems, such as  the question whether the Jones polynomial detects knots.  
Counterexamples to this conjecture exist in the virtual domain.  It is an
open problem whether some of these counterexamples are isotopic to  classical knots and links.

\section{The $L$-equivalence for Virtual Braids}

Just as classical knots and links can be represented by the closures of braids, so can  virtual
knots and links be represented by the closures of virtual braids \cite{SVKT,Ka,KL}.  A {\it virtual braid}
on $n$ strands is a braid on $n$ strands in the classical sense, which may also contain virtual
crossings. The closure of a virtual braid is formed by joining by simple
arcs the corresponding endpoints of the braid on its plane.  
It is easily seen that the classical Alexander Theorem  generalizes to
virtuals \cite{Ka,KL}. 

\begin{th}{ \ Every 
(oriented) virtual link  can be represented by a virtual braid, whose closure is
isotopic to the original link.
}
\end{th}
 
\bigbreak 

$$\vbox{\picill1.3inby2.3in(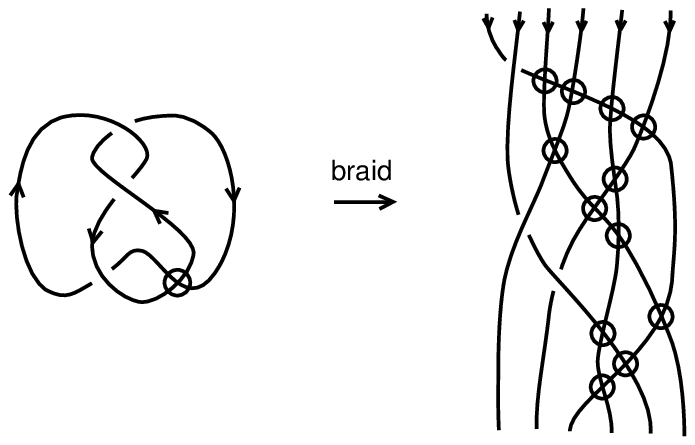)  }$$

\begin{center}
{ \bf Figure 4 --  An Example of Braiding} 
\end{center}

As in classical knot theory, the next consideration after the braiding is to characterize virtual
braids that induce, via closure, isotopic virtual links. In this section we define $L$--equivalence of virtual braids. 
For this purpose we
need to recall and generalize to the virtual setting the {\it $L$--move} between braids. The  $L$--move was
introduced in \cite{La,LR}, where it was used for proving the ``one--move Markov theorem" for classical
oriented links (cf.Theorem 2.3 in \cite{LR}), where it replaces the two well-known moves of the Markov
equivalence: the {\it stabilization} that introduces a crossing at the bottom right of a braid and
{\it conjugation} that conjugates a braid by a crossing. 

\begin{defn}{\rm \ 
A  {\it basic virtual $L$--move}  on a virtual braid, denoted  {\it  $L_v$--move}, 
consists in cutting an arc of the  braid open and pulling the upper cutpoint downward and
the lower  upward, so as to create a new pair of braid strands with corresponding endpoints
(on the vertical line of the cutpoint), and such that both strands  cross entirely  {\it virtually}
 with the rest of the braid. (In abstract illustrations this is indicated by placing virtual crossings on
the border of the braid box.) 

By a small braid isotopy that does not  change the relative positions of endpoints,
an $L_v$--move can be equivalently seen as introducing an in--box virtual crossing to a virtual braid,
which faces either the {\it right} or the {\it left} side of the braid. If we want to
emphasize the existence of the virtual crossing, we shall say {\it virtual $L_v$-move}, denoted {\it $vL_v$--move}. 
In  Figure 5 we give abstract illustrations. See also Figure 10 for a concrete example.}
\end{defn}
\bigbreak

$$\vbox{\picill4inby1.7in(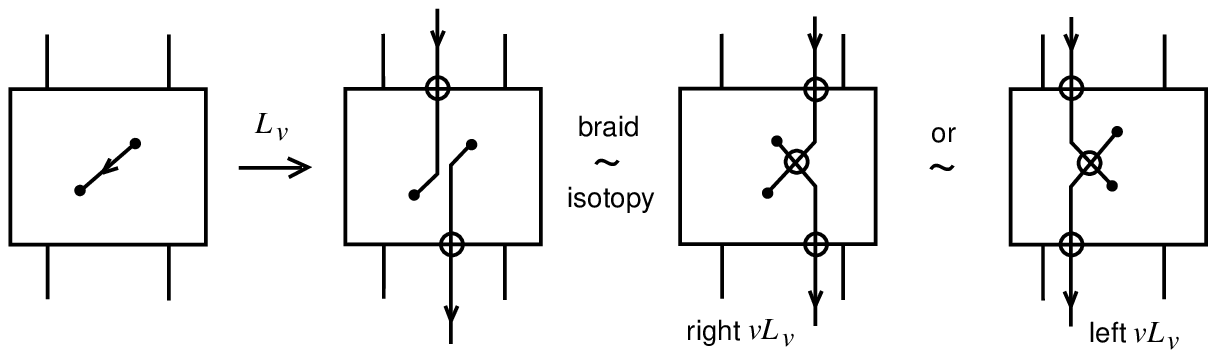)  }$$ 

\begin{center}
{ \bf Figure 5 -- A Basic Virtual $L$--move Without and With a Virtual Crossing} 
\end{center}

 Note  that in the closure of a $vL_v$--move the detoured loop contracts to a kink. This kink could also
be created by a real crossing, positive or negative. So we have:

\begin{defn}{\rm \ 
A  {\it real $L_v$--move}, abbreviated to  {\it  $+L_v$--move} or {\it  $-L_v$--move}, is a virtual 
$L$--move that introduces  a real in--box crossing  on a virtual braid, and it can face either the {\it
right} or the {\it left} side of the braid. See  Figure 6 for an illustration.}
\end{defn}

\bigbreak

$$\vbox{\picill4.2inby1.7in(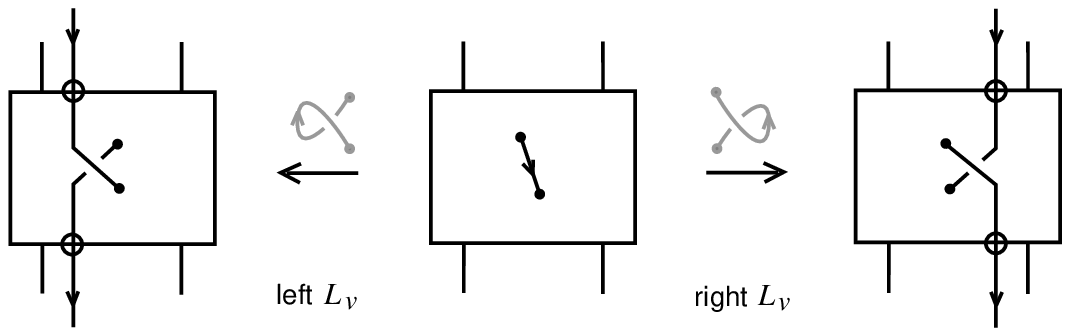)  }$$ 

\begin{center}
{ \bf Figure 6 --  Left and Right  Real $L_v$--moves} 
\end{center}

If the crossing of the kink is virtual, then, in the presence of the forbidden moves, there
is another possibility for a  move  on the braid level, which uses another arc of the braid, the `thread'. 

\begin{defn}{\rm \ 
A   {\it threaded virtual $L$--move} on a virtual braid is a virtual  $L$--move with a virtual
crossing in which, before pulling open the little up--arc of the kink, we perform a Reidemeister II
move with real crossings, using another arc of the braid, the {\it thread}. There are two possibilities: A 
{\it threaded over $L_v$--move} and a {\it threaded under $L_v$--move}, depending on whether we pull the kink 
over or under the thread,  both with the variants {\it right and left}.  See Figure 7 for abstract
illustrations.  }
\end{defn}

\bigbreak

$$\vbox{\picill4.2inby1.7in(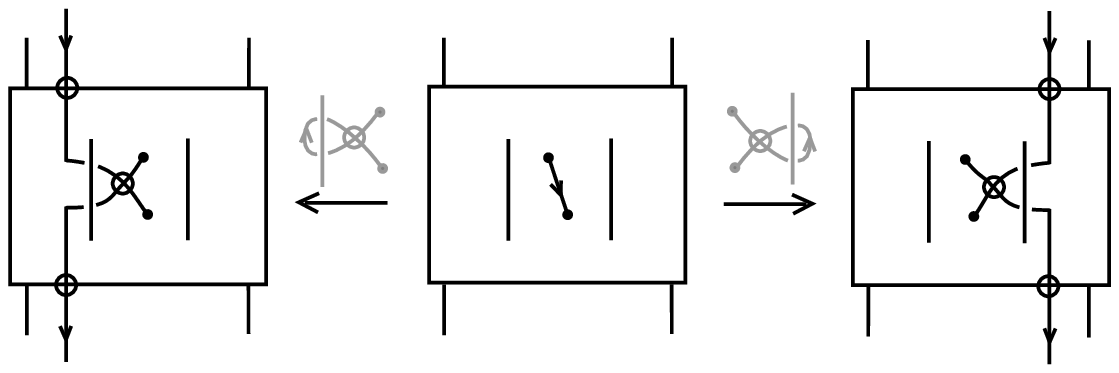)  }$$ 

\begin{center}
{ \bf Figure 7 -- Left and Right Threaded Under $L_v$--moves} 
\end{center}

Note that a threaded virtual $L$--move cannot be simplified in the braid. 
 If the crossing of the kink were real then, using a braid RIII move  
with the thread, the move would reduce  to a virtual $L$--move with a real crossing. Similarly, if the
forbidden moves were allowed, a threaded virtual $L$--move would reduce to a  $vL_v$--move.

\begin{rem}{\rm \  As with a braiding move, the effect of a  virtual $L$--move, basic, real or threaded, is to 
stretch an arc of the braid around the braid axis using the detour move after twisting it and possibly after
threading it. On the other hand, such a move between virtual braids gives rise to isotopic closures,
since the virtual $L$--moves shrink locally to kinks (grey diagrams in Figures 6 and 7). }
\end{rem}

 Conceivably, the threading of a virtual $L$--move could involve a
sequence of threads and Reidemeister II moves with over, under or virtual crossings, as Figure 8
suggests. The presence of the forbidden moves does not allow for simplifications on the braid level.  We show
in \cite{VL} that  the {\it multi--threaded $L_v$--moves}  follow from the simple threaded moves. 

\bigbreak

$$\vbox{\picill4inby1.7in(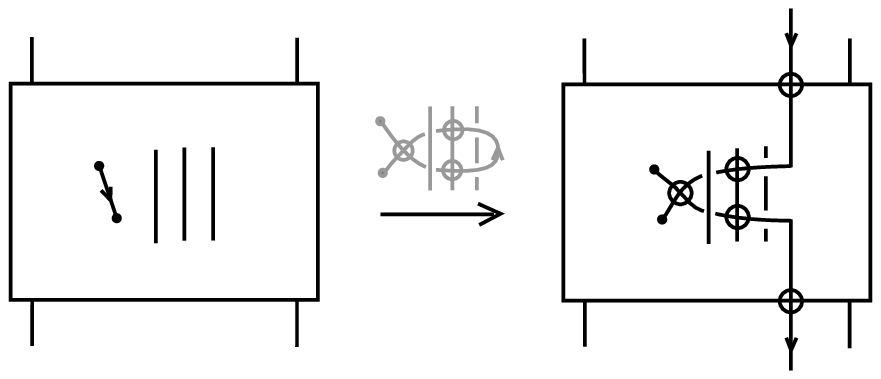)  }$$ 

\begin{center}
{ \bf Figure 8 -- A Right Multi-Threaded Virtual $L$--move} 
\end{center}

We finally introduce the notion of a classical $L$--move, adapted to our set--up. 

\begin{defn}{\rm \ A {\it classical $L_{over}$--move} resp.  {\it $L_{under}$--move } on a virtual braid
consists in  cutting an arc of the virtual braid open and pulling the two ends, so as to create a new
pair of braid strands, which run both {\it entirely over}  resp. {\it entirely under}  the rest of the
braid, and such that the closures of the virtual braids before and after the move are isotopic. See
Figure 9 for abstract illustrations.  A classical  $L$--move may also introduce an in--box
crossing, which may be positive, negative or virtual,  or it may even involve a thread. }
\end{defn}

In order that a classical $L$--move between virtual braids is {\it allowed}, in the sense that it gives rise to
isotopic virtual links upon closure, it is required that the virtual braid  has no virtual crossings on the
entire vertical zone either to the left or to the right of the new strands of the $L$--move. We then place the
axis of the braid perpendicularly to the plane, on the side with no virtual crossings. 
 We show in \cite{VL} that the allowed  $L$--moves can be expressed in terms of  virtual
$L$--moves and real conjugation. It was the classical $L$--moves that were introduced  in \cite{LR}, and they
replaced the two equivalence moves of the classical Markov theorem. Clearly, in the classical set--up these
moves are always allowed, while the presence of forbidden moves can preclude them in the virtual setting.

\bigbreak

$$\vbox{\picill4inby1.6in(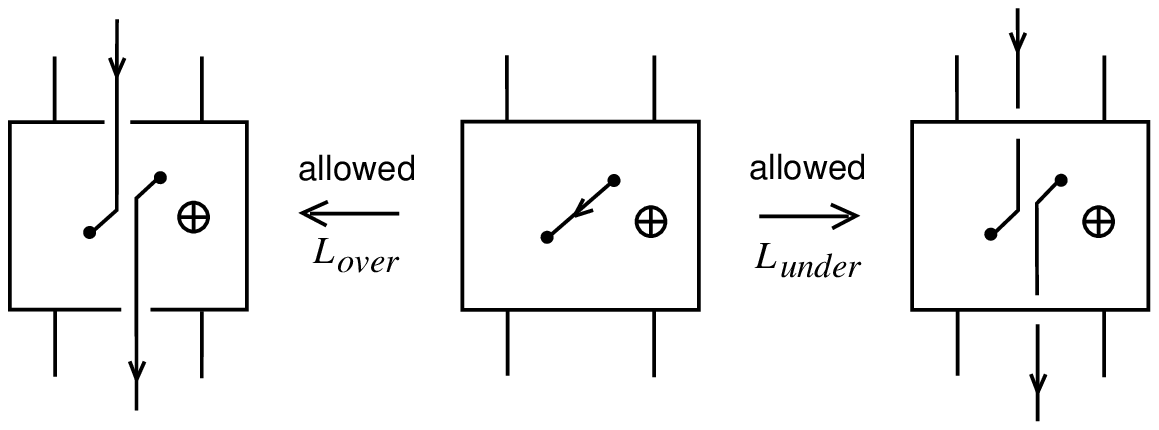)  }$$ 

\begin{center}
{ \bf Figure 9 -- The Allowed Classical  $L$--moves} 
\end{center}

 In Figure 10  we illustrate an example of various types of $L$--moves taking place at the same point
of a virtual  braid.

\bigbreak

\bigbreak 

$$\vbox{\picill5inby1.5in(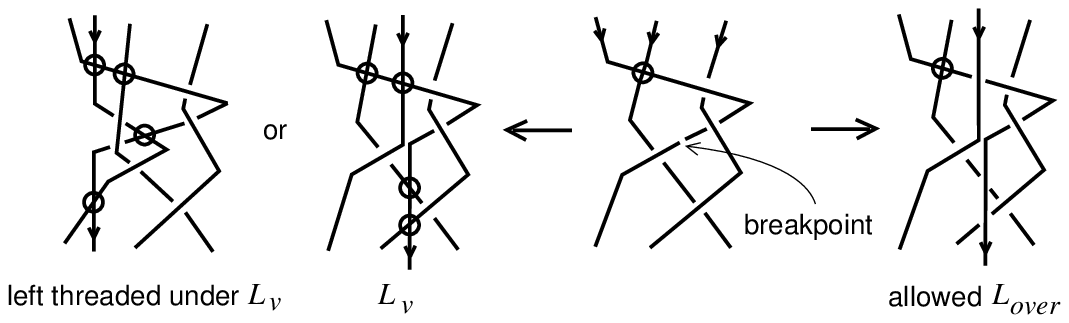)  }$$ 

\begin{center}
{ \bf Figure 10 -- A Concrete Example of Introducing $L$--moves} 
\end{center}

\section {The $L$--move Markov Theorem for Virtual Braids}

It is clear that different choices when applying the braiding algorithm  as well as local 
isotopy changes on the diagram level may result in  different virtual braids. 
In the proof of the following theorem we show that {\it real conjugation} (that is, conjugation by a real crossing)
and the $L_v$--moves with all variations (recall Definitions $1,2$ and $3$) capture and reflect on the braid
level all instances of isotopy between virtual links.

\begin{th}[$L$--move Markov Theorem for virtuals] Two oriented  virtual links are isotopic if and
only if any two corresponding virtual braids differ by virtual braid isotopy and a finite sequence of the
following moves or their inverses:
\begin{itemize}
\item[(i)]Real conjugation 
\item[(ii)]Right virtual $L_v$--moves  
\item[(iii)]Right real $L_v$--moves  
\item[(iv)]Right and left threaded $L_v$--moves.
\end{itemize} 
\end{th}

\noindent Moves  (ii), (iii), (iv) together  with their
inverses shall be called collectively {\it virtual $L$--moves}. Virtual $L$--moves together with virtual
braid isotopy generate an equivalence relation in the set of virtual braids, the {\it $L$--equivalence}.
\bigbreak

\section {Algebraic Markov Equivalence for Virtual Braids}

In this section  we reformulate and sharpen the statement of Theorem~4.1 by giving an equivalent list of local
algebraic moves in the virtual braid groups. More precisely, let  $VB_{n}$ denote the virtual braid group on
$n$ strands and let $\sigma_i, v_i$ be its generating classical and virtual crossings. The $\sigma_i$'s
satisfy the relations of the classical braid group and the $v_i$'s satisfy the relations of the permutation
group. The characteristic relation in $VB_{n}$  is the {\it mixed relation} relating both:
$ \begin{array}{cccl} 
v_{i} \sigma_{i+1} v_{i} & = & v_{i+1} \sigma_i v_{i+1}. &   \\ 
\end{array}$
\bigbreak

$$\vbox{\picill4inby6in(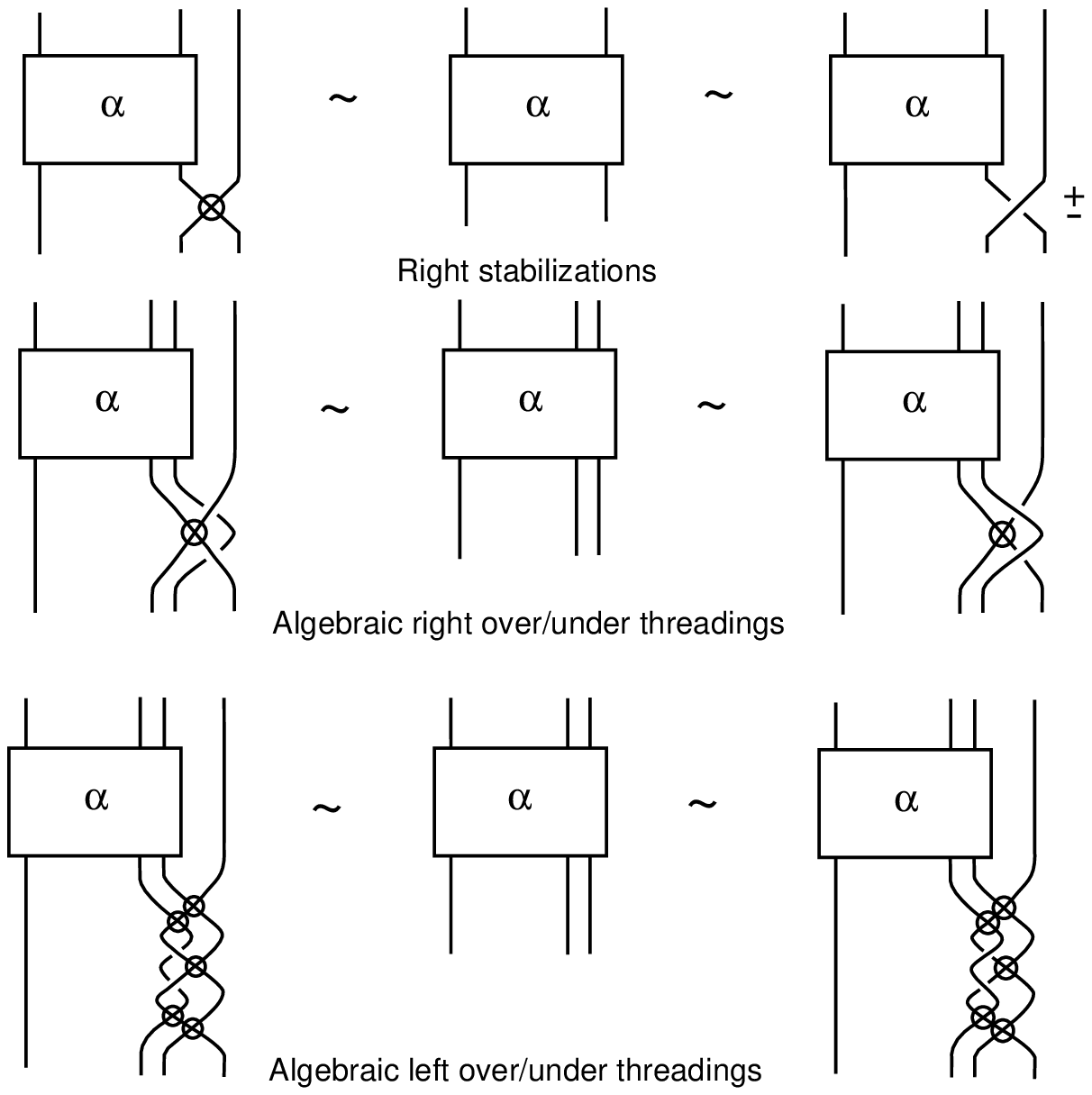)  }$$ 

\begin{center}
{ \bf Figure 11 -- The Moves (ii), (iii) and (iv) of Theorem 3} 
\end{center}

\noindent The group $VB_{n}$ embedds naturally into $VB_{n+1}$ by adding one identity strand at the right of
the braid. So, it makes sense to define $VB_{\infty} :=
\bigcup_{n=1}^{\infty} VB_{n}$,  the disjoint union of all virtual braid groups. We can now state our result.

\begin{th}[Algebraic Markov Theorem for virtuals] Two oriented  virtual links are isotopic if and
only if any two corresponding virtual braids differ by braid relations in $VB_{\infty}$ and a finite
sequence of the following moves or their inverses:

\begin{itemize}
\item[(i)]Virtual and real conjugation:    \ \ \ \ \ \ \ \ \ \ \ \  $ v_i \alpha v_i \sim \alpha \sim 
{\sigma_i}^{-1}\alpha \sigma_i $
\item[(ii)]Right virtual and real stabilization:  \ \ \ \ \  $\alpha v_n \sim \alpha
\sim \alpha \sigma_n^{\pm 1}$ 
\item[(iii)]Algebraic right over/under threading:  \ \  $\alpha \sim \alpha \sigma_n^{\pm 1} v_{n-1}
\sigma_n^{\mp 1} $
\item[(iv)]Algebraic left over/under threading: $$\alpha \sim  \alpha v_n v_{n-1} \sigma_{n-1}^{\mp
1} v_n \sigma_{n-1}^{\pm 1} v_{n-1} v_n $$
\end{itemize} 

\noindent where $\alpha,  v_i, \sigma_i \in VB_n$ and  $v_n, \sigma_n \in VB_{n+1}$ (see Figure 11).
\end{th}

\bigbreak

 \noindent {\bf Acknowledgments.} \ We are happy to mention that the paper of S. Kamada \cite{Ka}
has been for us a source of inspiration. It gives the first author great pleasure to acknowledge support from  NSF Grant
DMS-0245588.
\bigbreak

\end{document}